\documentclass{article}%
\usepackage{amsmath}
\usepackage{graphicx}%
\usepackage{amsfonts}%
\usepackage{amssymb}

\begin{document}

\title{Fourier sine and cosine transforms of irrational functions}
\author{Bernard J. Laurenzi}
\date{December 9, 2019\\
\ \ \ \ \ \ \ \ Department of Chemistry \\
\ \ \ \ \ \ \ The State University of New York at Albany,\\
Albany, New York, 12222}
\maketitle

\begin{abstract}
Fourier sine transforms containing irrational integrands are presented.
\ Explicit closed form expressions are shown to be related to Lommel functions
and in special cases to the Fresnel integrals.

\end{abstract}

\section{Fourier Transforms of some irrational integrands}

\ \ \ \ \ \ \ \ \ \ \ \ \ \ \ \ \ \ \ \ \ \ \ \ \ \ \ \ \ \ \ \ \ \ \ \ \ \ \ \ \ \ \ \ \ \ \ \ \ \ \ \ \ \ \ 

Recently integrals of the form%
\[
\int_{0}^{\infty}\frac{\sin(\zeta t)}{\sqrt{(t+a)}\sqrt{t+b}}dt,\hspace
{0.25in}\int_{0}^{\infty}\frac{\cos(\zeta t)}{\sqrt{(t+a)}\sqrt{t+b}}dt,
\]
have arisen in a semi-classical theory of atomic structure \cite{Englert}.
\ Since there appears to be little known about the exact values for these
integrals we present the forms shown below. \ 

Before investigation of the integrals shown above we begin with a family of
the less complex integrals i.e.
\begin{subequations}
\begin{align}
\mathbf{S}_{\alpha}(x,\zeta)  &  =\int_{0}^{\infty}\frac{\sin(\zeta
t)}{(t+x)^{\alpha+1/2}}dt,\label{eq1}\\
\mathbf{C}_{\alpha}(x,\zeta)  &  =\int_{0}^{\infty}\frac{\cos(\zeta
t)}{(t+x)^{\alpha+1/2}}dt, \label{eq2}%
\end{align}
for which closed form expressions for $\mathbf{S}_{0}(x,\zeta)$ and
$\mathbf{C}_{0}(x,\zeta)$ are known \cite{Laurenzi} i.e.
\end{subequations}
\begin{align*}
\mathbf{S}_{0}(x,\zeta)  &  =\sqrt{\frac{\pi}{2\zeta}}[\cos(\zeta
x)\{1-2S(\sqrt{\frac{2\zeta x}{\pi}})\}-\sin(\zeta x)\{1-2C(\sqrt{\frac{2\zeta
x}{\pi}})\}],\\
\mathbf{C}_{0}(x,\zeta)  &  =\sqrt{\frac{\pi}{2\zeta}}[\cos(\zeta
x)\{1-2C(\sqrt{\frac{2\zeta x}{\pi}})\}+\sin(\zeta x)\{1-2S(\sqrt{\frac{2\zeta
x}{\pi}})\}],
\end{align*}
and where $S(z)$ and $C(z)$ are the Fresnel integrals \cite{NBS}%
\begin{align*}
S(z)  &  =\int_{0}^{z}\sin(\frac{\pi}{2}t^{2})dt,\\
C(z)  &  =\int_{0}^{z}\cos(\frac{\pi}{2}t^{2})dt.
\end{align*}
We note that the $\mathbf{S}_{0}(x,\zeta)$ and $\mathbf{C}_{0}(x,\zeta)$
integrals are also related to the Lommel functions \cite{Lommel} of the second
kind $S_{\mu,1/2}(z)$. \ In the case where $\alpha$ is any real number,
generalization of the integrals in (1a,1b) can be written in terms of the
Lommel functions. cf. appendix A.

As will be seen below, closed form expressions for the integrals with integer
values of $\alpha>0$\ have been obtained from $\mathbf{S}_{0}(x,\zeta)$ and
$\mathbf{C}_{0}(x,\zeta)$, this family of integrals being the main result of
the present work.\ \ We note that the integrals $\mathbf{S}_{\alpha}(x,\zeta)$
and $\mathbf{C}_{\alpha}(x,\zeta)$ scale as
\begin{align*}
\mathbf{S}_{\alpha}(x,\zeta)  &  =\zeta^{\alpha-1/2}\,\mathbf{S}_{\alpha
}(\zeta x,1),\\
\mathbf{C}_{\alpha}(\zeta x,\zeta)  &  =\zeta^{\alpha-1/2}\mathbf{C}_{\alpha
}(\zeta x,1).
\end{align*}
Differentiating the scaled integrals with respect to $\zeta x$ we have%
\begin{align*}
\frac{d\,\mathbf{S}_{\alpha}(\zeta x,1)}{d(\zeta x)}  &  =-(\alpha
+1/2)\,\mathbf{S}_{\alpha+1}(\zeta x,1),\\
\frac{d\mathbf{C}_{\alpha}(\zeta x,1)}{d(\zeta x)}  &  =-(\alpha
+1/2)\mathbf{C}_{\alpha+1}(\zeta x,1),
\end{align*}
and
\begin{subequations}
\begin{align}
\frac{d^{\,2}\,\mathbf{S}_{\alpha}(\zeta x,1)}{d(\zeta x)^{2}}  &
=(\alpha+1/2)(\alpha+3/2)\mathbf{S}_{\alpha+2}(\zeta x,1),\label{2q2}\\
\,\,\frac{d^{\,2}\,\mathbf{C}_{\alpha}(\zeta x,1)}{d(\zeta x)^{2}}  &
=(\alpha+1/2)(\alpha+3/2)\mathbf{C}_{\alpha+2}(\zeta x,1). \label{eq3}%
\end{align}
Integration by parts of the integrals $\mathbf{S}_{\alpha}(\zeta x,1)$ and
$\mathbf{C}_{\alpha}(x,1)$ give relations which interrelate these functions
i.e.
\end{subequations}
\begin{align*}
\mathbf{S}_{\alpha}(\zeta x,1)  &  =\frac{1}{(\zeta x)^{\alpha+1/2}}%
-(\alpha+1/2)\mathbf{C}_{\alpha+1}(\zeta x,1),\\
\mathbf{C}_{\alpha}(\zeta x,1)  &  =(\alpha+1/2)\mathbf{S}_{\alpha+1}(\zeta
x,1).
\end{align*}
Eliminating $C_{\alpha}$ from the equations above we get the difference
equation
\begin{equation}
(\alpha+1/2)(\alpha+3/2)\mathbf{S}_{\alpha+2}(\varsigma x,1)+\mathbf{S}%
_{\alpha}(\zeta x,1)=\frac{1}{(\zeta x)^{\alpha+1/2}}, \label{eq4}%
\end{equation}
and similarly for its counterpart $\mathbf{C}_{\alpha}$ we have
\[
(\alpha+1/2)(\alpha+3/2)\mathbf{C}_{\alpha+2}(\zeta x,1)+\mathbf{C}_{\alpha
}(\zeta x,1)=\frac{(\alpha+1/2)}{(\zeta x)^{\alpha+3/2}}.
\]
Using there results together with equations (2a) and (2b) we get the
differential equations%
\begin{align*}
\frac{d^{\,2}\,\mathbf{S}_{\alpha}(\zeta x,1)}{d(\zeta x)^{2}}+\mathbf{S}%
_{\alpha}(\zeta x,1)  &  =\frac{1}{(\zeta x)^{\alpha+1/2}},\\
\,\,\frac{d^{\,2}\,\mathbf{C}_{\alpha}(\zeta x,1)}{d(\zeta x)^{2}}%
+\mathbf{C}_{\alpha}(\zeta x,1)  &  =\frac{(\alpha+1/2)}{(\zeta x)^{\alpha
+3/2}}.
\end{align*}

\subsection{The integrals $\mathbf{S}_{\alpha}(\zeta x,1)$ and $\mathbf{C}%
_{\alpha}(\zeta x,1)$ in the case $\alpha=2n$}

\ \ \ \ \ \ \ \ \ \ \ \ \ \ \ \ \ \ \ \ \ \ \ \ \ \ \ \ \ \ \ \ \ \ \ \ \ \ \ \ \ \ \ \ \ \ \ \ \ \ \ \ \ \ \ 

\bigskip In the case where $\alpha$ is an even integer we write with the help
of the recurrence relation for $\mathbf{S}_{\alpha}(\zeta x,1)$
\[
\mathbf{S}_{2n}(\zeta x,1)=F_{2n}(\zeta x)+c_{2n}[\cos(\zeta x)\{1-2S(\sqrt
{\frac{2\zeta x}{\pi}})\}-\sin(\zeta x)\{1-2C(\sqrt{\frac{2\zeta x}{\pi}%
})\}],
\]
where $F_{2n}(\zeta x)$ are functions of $\zeta x$ and the quantities $c_{2n}$
are constants. \ Using (3) we have
\[
(4n+1)(4n+3)F_{2n+2}(\zeta x)+4F_{2n}(\zeta x)=4/(\zeta x)^{2n+1/2},
\]
and
\[
(4n+1)(4n+3)c_{2n+2}+4c_{2n}=0.
\]
The solutions to these \textit{first-order} difference equations with initial
conditions $c_{0}=\sqrt{\pi/2}$ and $F_{0}(\zeta x)=0$ are
\[
F_{2n}(\zeta x)=(-1)^{n+1}\sum_{k=0}^{n-1}\frac{(-1)^{k}}{(\zeta x)^{2k+1/2}%
}\frac{\Gamma(2k+1/2)}{\Gamma(2n+1/2)},
\]
and%
\[
c_{2n}=\frac{(-1)^{n}\sqrt{2}\pi}{2\Gamma(2n+1/2)}.
\]

Similarly in the case for $\mathbf{C}_{2n}(x,1)$ we have with $\widehat{c}%
_{0}=\sqrt{\pi/2},$ and $\widehat{F}_{0}(\zeta x)=0$ and the form%
\[
\mathbf{C}_{2n}(\zeta x,1)=\widehat{F}_{2n}(\zeta x)+\widehat{c}_{2n}%
[\cos(\zeta x)\{1-2C(\sqrt{\frac{2\zeta x}{\pi}})\}+\sin(\zeta x)\{1-2S(\sqrt
{\frac{2\zeta x}{\pi}})\}],
\]
from which follows the difference equations%
\[
(4n+1)(4n+3)\widehat{F}_{2n+2}(\zeta x)+4\widehat{F}_{2n}(\zeta
x)=(8n+2)/(\zeta x)^{2n+1/2},
\]%
\[
(4n+1)(4n+3)\widehat{c}_{2n+2}+4\widehat{c}_{2n}=0.
\]
Here the corresponding solutions are%
\[
\widehat{F}_{2n}(\zeta x)=(-1)^{n+1}\sum_{k=0}^{n-1}\frac{(-1)^{k}}{(\zeta
x)^{2k+3/2}}\frac{\Gamma(2k+3/2)}{\Gamma(2n+1/2)},
\]
and%
\[
\widehat{c}_{2n}=\frac{(-1)^{n}\sqrt{2}\pi}{2\Gamma(2n+1/2)}.
\]

\subsection{The integrals $\mathbf{S}_{\alpha}(\zeta x,1)$ and $\mathbf{C}%
_{\alpha}(\zeta x,1)$ in the case $\alpha=2n+1$}

\ \ \ \ \ \ \ \ \ \ \ \ \ \ \ \ \ \ \ \ \ \ \ \ \ \ \ \ \ \ \ \ \ \ \ \ \ \ \ \ \ \ \ \ \ \ \ \ \ \ \ \ \ \ \ \ \ \ \ \ \ \ \ \ \ \ 

\bigskip

In the case where $\alpha$ is an odd integer we have
\[
\mathbf{S}_{2n+1}(\zeta x,1)=F_{2n+1}(\zeta x)+c_{2n+1}[\cos(\zeta
x)\{1-2C(\sqrt{\frac{2\zeta x}{\pi}})\}+\sin(\zeta x)\{1-2S(\sqrt{\frac{2\zeta
x}{\pi}})\}],
\]
with the difference equations%
\[
(4n+3)(4n+5)F_{2n+3}(\zeta x)+4F_{2n+1}(\zeta x)=4/(\zeta x)^{2n+3/2},
\]
and%
\[
(4n+3)(4n+5)c_{2n+3}+4c_{2n+1}=0.
\]
The latter solutions are with $c_{1}=\sqrt{2\pi}$ and $F_{1}(\zeta x)=0$
\[
F_{2n+1}(\zeta x)=(-1)^{n}\sum_{k=0}^{n-1}\frac{(-1)^{k}}{(\zeta x)^{2k+3/2}%
}\frac{\Gamma(2k+3/2)}{\Gamma(2n+3/2)},
\]
and%
\[
c_{2n+1}=\frac{(-1)^{n}\sqrt{2}\pi}{2\Gamma(2n+3/2)}.
\]
In the case of $\mathbf{C}_{2n+1}(\zeta x,1)$ we have%
\[
\mathbf{C}_{2n+1}(\zeta x,1)=\widehat{F}_{2n+1}(\zeta x)+\widehat{c}%
_{2n+1}[\cos(\zeta x)\{1-2S(\sqrt{\frac{2\zeta x}{\pi}})\}-\sin(\zeta
x)\{1-2C(\sqrt{\frac{2\zeta x}{\pi}})\}],
\]
with difference equations%
\[
(4n+3)(4n+5)\widehat{F}_{2n+3}(\zeta x)+4\widehat{F}_{2n+1}(\zeta
x)=(8n+6)/(\zeta x)^{2n+5/2},
\]%
\[
(4n+3)(4n+5)\widehat{c}_{2n+3}+4\widehat{c}_{2n+1}=0,
\]
the solutions with $\widehat{c}_{1}=-\sqrt{2\pi}$ and $\widehat{F}_{1}(\zeta
x)=2/\sqrt{\zeta x}$ \ being%
\[
\widehat{F}_{2n+1}(\zeta x)=\frac{(-1)^{n+1}\sqrt{\pi}}{\Gamma(2n+3/2)\sqrt
{\zeta x}}+(-1)^{n+1}\sum_{k=0}^{n-1}\frac{(-1)^{k}}{(\zeta x)^{2k+5/2}%
}\frac{\Gamma(2k+5/2)}{\Gamma(2n+3/2)},
\]
and%

\[
\widehat{c}_{2n+1}=\frac{(-1)^{n+1}\sqrt{2}\pi}{2\Gamma(2n+3/2)}.
\]

\subsection{The irrational integrals $\int_{0}^{\infty}\tfrac{\sin(\zeta
t)}{\sqrt{t+a}\sqrt{t+b}}dt$ and $\int_{0}^{\infty}\frac{\cos(\zeta t)}%
{\sqrt{(t+a)}\sqrt{t+b}}dt.$}

\ \ \ \ \ \ \ \ \ \ \ \ \ \ \ \ \ \ \ \ \ \ \ \ \ \ \ \ \ \ \ \ \ \ \ \ \ \ \ \ \ \ \ \ \ \ \ \ \ \ \ \ \ \ \ \ \ \ \ \ \ \ \ \ \ \ 

\bigskip We now take up the integrals%
\begin{equation}
\int_{0}^{\infty}\frac{\sin(\zeta t)}{\sqrt{t+a}\sqrt{t+b}}dt,\hspace
{0.25in}\int_{0}^{\infty}\frac{\cos(\zeta t)}{\sqrt{(t+a)}\sqrt{t+b}}dt.
\label{eq5}%
\end{equation}
With a change of variable we have
\begin{align*}
\int_{0}^{\infty}\frac{\sin(\zeta t)}{\sqrt{t+a}\sqrt{t+b}}dt  &
=2\cos(a\zeta)\int_{\sqrt{a}}^{\infty}\frac{\sin(\zeta s^{2})}{\sqrt
{s^{2}+(b-a)}}ds-2\sin(a\zeta)\int_{\sqrt{a}}^{\infty}\frac{\cos(\zeta s^{2}%
)}{\sqrt{s^{2}+(b-a)}}ds,\\
\int_{0}^{\infty}\frac{\cos(\zeta t)}{\sqrt{t+a}\sqrt{t+b}}dt  &
=2\cos(a\zeta)\int_{\sqrt{a}}^{\infty}\frac{\cos(\zeta s^{2})}{\sqrt
{s^{2}+(b-a)}}ds+2\sin(a\zeta)\int_{\sqrt{a}}^{\infty}\frac{\sin(\zeta s^{2}%
)}{\sqrt{s^{2}+(b-a)}}ds.
\end{align*}
Scaling those integrals we obtain with $b>a$%
\begin{equation}
\int_{0}^{\infty}\frac{\sin(\zeta t)}{\sqrt{t+a}\sqrt{t+b}}dt=2\cos
(a\zeta)\int_{\gamma}^{\infty}\frac{\sin(\zeta\lbrack b-a]z^{2})}{\sqrt
{z^{2}+1}}dz-2\sin(a\zeta)\int_{\gamma}^{\infty}\frac{\cos(\zeta\lbrack
b-a]z^{2})}{\sqrt{z^{2}+1}}dz \label{eq6}%
\end{equation}
and%
\begin{equation}
\int_{0}^{\infty}\frac{\cos(\zeta t)}{\sqrt{t+a}\sqrt{t+b}}dt=2\cos
(a\zeta)\int_{\gamma}^{\infty}\frac{\cos(\zeta\lbrack b-a]z^{2})}{\sqrt
{z^{2}+1}}dz+2\sin(a\zeta)\int_{\gamma}^{\infty}\frac{\sin(\zeta\lbrack
b-a]z^{2})}{\sqrt{z^{2}+1}}dz, \label{eq7}%
\end{equation}
where $\gamma=\sqrt{\frac{a}{b-a}}.$

The integrals (5,6) with \ infinite range are well known and are given by
\[
\int_{0}^{\infty}\frac{\sin(\zeta\lbrack b-a]z^{2})}{\sqrt{z^{2}+1}%
}dz=\frac{\pi}{4}\{\sin(\frac{\zeta\lbrack b-a]}{2})Y_{0}(\frac{\zeta\lbrack
b-a]}{2})+\cos(\frac{\zeta\lbrack b-a]}{2})J_{0}(\frac{\zeta\lbrack b-a]}%
{2})\},
\]
and%
\[
\int_{0}^{\infty}\frac{\cos(\zeta\lbrack b-a]z^{2})}{\sqrt{z^{2}+1}%
}dz=\frac{\pi}{4}\{\sin(\frac{\zeta\lbrack b-a]}{2})J_{0}(\frac{\zeta\lbrack
b-a]}{2})-\cos(\frac{\zeta\lbrack b-a]}{2})Y_{0}(\frac{\zeta\lbrack b-a]}%
{2})\},
\]
where $J_{0}(z)$ and $Y_{0}(z)$ are Bessel functions \cite{NBS2} of the first
and second kind with degree zero. \ The corresponding integrals with finite
range i.e.
\[
\int_{0}^{\sqrt{\frac{a}{b-a}}}\frac{\sin(\zeta\lbrack b-a]z^{2})}{\sqrt
{z^{2}+1}}dz,\hspace{0.25in}\int_{0}^{\sqrt{\frac{a}{b-a}}}\frac{\cos
(\zeta\lbrack b-a]z^{2})}{\sqrt{z^{2}+1}}dz,
\]
can be obtained as rapidly converging infinite series by expanding the $\sin$
and $\cos$ functions in Taylor series and integrating term by term. \ We get
with $c=\zeta\lbrack b-a]$ and $\gamma=\sqrt{\frac{a}{b-a}}\leq1,$%
\[
\int_{0}^{\sqrt{\frac{a}{b-a}}}\frac{\sin(\zeta\lbrack b-a]z^{2})}{\sqrt
{z^{2}+1}}dz=c\gamma^{3}\sum_{k=0}^{\infty}\frac{(-c^{2}\gamma^{4})^{k}%
}{(2k+1)!(4k+3)}\,_{2}F_{1}(1/2,2k+3/2;2k+5/2;-\gamma^{2}),
\]
and%
\[
\int_{0}^{\sqrt{\frac{a}{b-a}}}\frac{\cos(\zeta\lbrack b-a]z^{2})}{\sqrt
{z^{2}+1}}dz=\gamma\sum_{k=0}^{\infty}\frac{(-c^{2}\gamma^{4})^{k}%
}{(2k)!(4k+1)}\,_{2}F_{1}(1/2,2k+1/2;2k+3/2;-\gamma^{2}).
\]
where $_{2}F_{1}$ are hypergeometric functions \cite{NBS3}. \ 

For $\gamma\leq1$ and keeping only the leading terms in the $\gamma$ series,
expansion of the hypergeometric functions shown above produces sums which
give
\begin{align*}
\int_{0}^{\sqrt{\frac{a}{b-a}}}\frac{\sin(\zeta\lbrack b-a]z^{2})}{\sqrt
{z^{2}+1}}dz  &  \thickapprox\frac{\gamma}{4c}\cos(c\gamma^{2})+\sqrt
{\frac{\pi}{2c}}\{S(\gamma\sqrt{\frac{2c}{\pi}})-\frac{1}{4c}C(\gamma
\sqrt{\frac{2c}{\pi}})\},\\
\int_{0}^{\sqrt{\frac{a}{b-a}}}\frac{\cos(\zeta\lbrack b-a]z^{2})}{\sqrt
{z^{2}+1}}dz  &  \thickapprox-\frac{\gamma}{c}\sin(c\gamma^{2})+\sqrt
{\frac{\pi}{2c}}\{\frac{1}{4c}S(\gamma\sqrt{\frac{2c}{\pi}})+C(\gamma
\sqrt{\frac{2c}{\pi}})\}.
\end{align*}
where the Fresnel integrals appear once again. \ Differentiating the integrals
in (5,6) with respect to either $a$ or $b$ produces higher-order integrals of
this kind. \ 

\subsection{Integrals $\int_{0}^{\infty}\frac{\sin(\zeta t)}{\sqrt{t+a}%
\sqrt{t+b}\sqrt{t+c}}dt$ and $\int_{0}^{\infty}\frac{\cos(\zeta t)}%
{\sqrt{(t+a)}\sqrt{t+b}\sqrt{t+c}}dt$}

\ \ \ \ \ \ \ \ \ \ \ \ \ \ \ \ \ \ \ \ \ \ \ \ \ \ \ \ \ \ \ \ \ \ \ \ \ \ \ \ \ \ \ \ \ \ \ \ \ \ \ \ \ \ \ \ \ \ \ \ \ \ \ \ \ \ \bigskip

Integrals of the form%
\[
\int_{0}^{\infty}\frac{\sin(\zeta t)}{\sqrt{t+a}\sqrt{t+b}\sqrt{t+c}%
}dt,\hspace{0.25in}\int_{0}^{\infty}\frac{\cos(\zeta t)}{\sqrt{(t+a)}%
\sqrt{t+b}\sqrt{t+c}}dt,
\]
which can be regarded as generalizations of the Carlson \cite{Carlson}
integrals though of great interest do not appear to yield to the methods used
above. \ However, in the case where two of the constants are equal i.e.
\[
\int_{0}^{\infty}\frac{\sin(\zeta t)}{\sqrt{t+a}\,(t+b)}dt,\hspace{0.25in}%
\int_{0}^{\infty}\frac{\cos(\zeta t)}{\sqrt{t+a}\,(t+b)}dt.
\]
progress can be made. \ By means similar to those used above we have
\begin{align*}
\int_{0}^{\infty}\frac{\sin(\zeta t)}{\sqrt{t+a}(t+b)}dt  &  =2\cos
(a\zeta)\int_{\sqrt{a}}^{\infty}\frac{\sin(\zeta s^{2})ds}{(s^{2}%
+(b-a))}-2\sin(a\zeta)\int_{\sqrt{a}}^{\infty}\frac{\cos(\zeta s^{2}%
)ds}{(s^{2}+(b-a))},\\
\int_{0}^{\infty}\frac{\cos(\zeta t)}{\sqrt{t+a}(t+b)}dt  &  =2\cos
(a\zeta)\int_{\sqrt{a}}^{\infty}\frac{\cos(\zeta s^{2})ds}{(s^{2}%
+(b-a))}+2\sin(a\zeta)\int_{\sqrt{a}}^{\infty}\frac{\sin(\zeta s^{2}%
)ds}{(s^{2}+(b-a))}.
\end{align*}
With scaling the integrals become
\begin{align*}
\int_{0}^{\infty}\frac{\sin(\zeta t)}{\sqrt{t+a}(t+b)}dt  &  =\frac{2\cos
(a\zeta)}{\sqrt{b-a}}\int_{\gamma}^{\infty}\frac{\sin(\zeta\lbrack
b-a]x^{2})dx}{(x^{2}+1)}-\frac{2\sin(a\zeta)}{\sqrt{b-a}}\int_{\gamma}%
^{\infty}\frac{\cos(\zeta\lbrack b-a]x^{2})dx}{(x^{2}+1)},\\
\int_{0}^{\infty}\frac{\cos(\zeta t)}{\sqrt{t+a}(t+b)}dt  &  =\frac{2\cos
(a\zeta)}{\sqrt{b-a}}\int_{\gamma}^{\infty}\frac{\cos(\zeta\lbrack
b-a]x^{2})dx}{(x^{2}+1)}+\frac{2\sin(a\zeta)}{\sqrt{(b-a)}}\int_{\gamma
}^{\infty}\frac{\sin(\zeta\lbrack b-a]x^{2})dx}{(x^{2}+1)}.
\end{align*}
The\ corresponding infinite range integrals are given by
\begin{align*}
\int_{0}^{\infty}\frac{\sin(\zeta\lbrack b-a]x^{2})dx}{(x^{2}+1)}  &
=\frac{\pi}{2}\{\sin(c)[S(\sqrt{\frac{2c}{\pi}})+C(\sqrt{\frac{2c}{\pi}%
})-1]-\cos(c)[S(\sqrt{\frac{2c}{\pi}})-C(\sqrt{\frac{2c}{\pi}})]\},\\
\int_{0}^{\infty}\frac{\cos(\zeta\lbrack b-a]x^{2})dx}{(x^{2}+1)}  &
=\frac{\pi}{2}\{\cos(c)[S(\sqrt{\frac{2c}{\pi}}+C(\sqrt{\frac{2c}{\pi}%
})+1]+\sin(c)[S(\sqrt{\frac{2c}{\pi}})-C(\sqrt{\frac{2c}{\pi}})]\}+\sqrt
{\frac{2\pi}{c}}.
\end{align*}
and the finite integrals as expressed as infinite series are
\begin{align*}
\int_{0}^{\sqrt{\frac{a}{b-a}}}\frac{\sin(\zeta\lbrack b-a]x^{2})dx}%
{(x^{2}+1)}  &  =c\gamma\sum_{k=0}^{\infty}\frac{(-c^{2}\gamma^{4})^{k}%
}{(2k+1)!(4k+1)}\{1-\,_{2}F_{1}(1,2k+1/2;2k+3/2;-\gamma^{2})\},\\
\int_{0}^{\sqrt{\frac{a}{b-a}}}\frac{\cos(\zeta\lbrack b-a]x^{2})dx}%
{(x^{2}+1)}  &  =\gamma\sum_{k=0}^{\infty}\frac{(-c^{2}\gamma^{4})^{k}%
}{(2k)!(4k+1)}\,_{2}F_{1}(1,2k+1/2;2k+3/2;-\gamma^{2}).
\end{align*}
For $\gamma\leq1$ these integrals become%
\begin{align*}
\int_{0}^{\sqrt{\frac{a}{b-a}}}\frac{\sin(\zeta\lbrack b-a]x^{2})dx}%
{(x^{2}+1)}  &  \thickapprox\frac{\gamma}{2c}\cos(c\gamma^{2})+\sqrt
{\frac{\pi}{2c}}\{S(\gamma\sqrt{\frac{2c}{\pi}})-\frac{1}{2c}C(\gamma
\sqrt{\frac{2c}{\pi}})\},\\
\int_{0}^{\sqrt{\frac{a}{b-a}}}\frac{\cos(\zeta\lbrack b-a]x^{2})dx}%
{(x^{2}+1)}  &  \thickapprox-\frac{\gamma}{2c}\sin(c\gamma^{2})+\sqrt
{\frac{\pi}{2c}}\{\frac{1}{2c}S(\gamma\sqrt{\frac{2c}{\pi}})+C(\gamma
\sqrt{\frac{2c}{\pi}})\},
\end{align*}

\begin{center}
\appendix          Appendix A
\end{center}

Integrals of the form shown in (1a) and (1b) with exponent $2n+1/m$ and
$2n+1+1/m$ \ with integer $n$ and $m$ can be written (cf. Maple) in terms of
the Lommel functions of the second kind $S_{\mu,1/2}(x)$ i.e.%

\begin{subequations}
\begin{align*}
\int_{0}^{\infty}\frac{\sin(\zeta t)}{(t+x)^{2n+1/m}}dt  &  =\zeta
^{2n+1/m-1}\sqrt{\zeta x}S_{-(2n+1/m-1/2),1/2}(\zeta x),\\
\int_{0}^{\infty}\frac{\cos(\zeta t)}{(t+x)^{2n+1/m}}dt  &  =\frac{\zeta
^{2n+1/m-1}}{(2n+1/m-1)}\left[  \frac{1}{(\zeta x)^{2n+1/m-1}}-\sqrt{\zeta
x}S_{-(2n+1/m-3/2),1/2}(\zeta x)\right]  ,\\
\int_{0}^{\infty}\frac{\sin(\zeta t)}{(t+x)^{2n+1+1/m}}dt  &  =\frac{\zeta
^{2n+1/m}}{(2n+1/m-1)(2n+1/m)}\left\{  \frac{1}{(\zeta x)^{2n+1/m-1}}%
-\sqrt{\zeta x}S_{-(2n+1/m-3/2),1/2}(\zeta x)]\right\}  ,\\
\int_{0}^{\infty}\frac{\cos(\zeta t)}{(t+x)^{2n+1+1/m}}dt  &  =\frac{\zeta
^{2n+1/m}}{(2n+1/m)}\left\{  \frac{1}{(\zeta x)^{2n+1/m}}-\sqrt{\zeta
x}S_{-(2n+1/m-1/2),1/2}(\zeta x)]\right\}  .
\end{align*}
Using the Lommel relation \cite{rel}
\end{subequations}
\[
z^{\mu+3/2}-\sqrt{z}S_{\mu+2,1/2}(z)=[(\mu+1)^{2}-1/4]\sqrt{z}S_{\mu,1/2}(z),
\]
three of the integrals shown above reduce to the forms below and we have in
summary%
\begin{align*}
\int_{0}^{\infty}\frac{\sin(\zeta t)}{(t+x)^{2n+1/m}}dt  &  =\zeta
^{2n+1/m-1}\sqrt{\zeta x}S_{-(2n+1/m-1/2),1/2}(\zeta x),\\
\int_{0}^{\infty}\frac{\cos(\zeta t)}{(t+x)^{2n+1/m}}dt  &  =\zeta
^{2n+1/m-1}(2n+1/m)\sqrt{\zeta x}S_{-(2n+1/m+1/2),1/2}(\zeta x),\\
\int_{0}^{\infty}\frac{\sin(\zeta t)}{(t+x)^{2n+1+1/m}}dt  &  =\zeta
^{2n+1/m}\sqrt{\zeta x}S_{-(2n+1/m+1/2),1/2}(\zeta x),\\
\int_{0}^{\infty}\frac{\cos(\zeta t)}{(t+x)^{2n+1+1/m}}dt  &  =\zeta
^{2n+1/m}(2n+1+1/m)\sqrt{\zeta x}S_{-(2n+1/m+3/2),1/2}(\zeta x).
\end{align*}
These integrals do not appear to have occurred previously in the literature
and can be seen to be integral representations of the Lommel function
$S_{\mu,1/2}$. \ In the case where is $n=0$ \ and $m=2$ as seen above, the
Lommel functions were related to the the Fresnel integrals, here we see that
this is generally the case. \ 

As given by Maple, one has the general relation%
\begin{equation}
\int_{0}^{\infty}\frac{\sin(t)}{(t+x)^{\alpha}}dt=\sqrt{x}S_{1/2-\alpha
,1/2}(x). \label{eq8}%
\end{equation}
The derivative of that integral with respect to $\alpha$ gives a new
logarithmic integral i.e.%
\[
\sqrt{x}\frac{\partial S_{1/2-\alpha,1/2}(x)}{\partial\alpha}=\frac{\partial
}{\partial\alpha}\int_{0}^{\infty}\frac{\sin(t)}{(t+x)^{\alpha}}dt=-\int
_{0}^{\infty}\frac{\ln(t+x)\sin(t)}{(t+x)^{\alpha}}dt.
\]
Direct integration of (8) gives%
\[
\int_{0}^{\infty}\frac{\sin(t)}{(t+x)^{\alpha}}dt=\tfrac{1}{2}\left\{
\exp(-\tfrac{i}{2}[\pi\alpha+2x])\Gamma(-\alpha,-ix)+\exp(\tfrac{i}{2}%
[\pi\alpha+2x])\Gamma(-\alpha,ix)\right\}  ,
\]
where $\Gamma(a,z)$ is the incomplete gamma function \cite{inc}
\ Differentiation of the integral above with respect to $\alpha$ then gives
the relation where $\alpha$ has been set equal to $1/2$ as%
\[
\left[  \frac{\partial}{\partial\alpha}\int_{0}^{\infty}\frac{\sin
(t)}{(t+x)^{\alpha}}dt\right]  _{\alpha=1/2}=
\]%
\begin{align*}
&  -\sqrt{x}\ln(x)S_{0,1/2}(x)-\tfrac{1}{2}\pi^{3/2}\sin(x+\pi/4)+\pi
^{1/2}[\gamma+\ln(4x)]\cos(x+\pi/4)\\
&  +4\sqrt{x}\left\{  \sin(x)\operatorname{Re}[_{2}F_{2}%
(1/2,1/2;3/2,3/2;ix)]-\cos(x)\operatorname{Im}[_{2}F_{2}%
(1/2,1/2;3/2,3/2;ix)]\right\}  ,
\end{align*}
where $_{2}F_{2}$ is a generalized hypergeometric function, \cite{hyper}. From
which relation the derivative of the Lommel function with respect to its first
parameter and the logarithmic integral above can be computed.

We note that it is also possible to write the integrals appearing above as
\begin{align*}
\int_{0}^{\infty}\frac{\sin(\zeta t)}{(t+x)^{2n+1/m}}dt  &  =\zeta
^{2n+1/m-1}\left[  \cos(\zeta x)\,si(-2n-1/m+1,\zeta x)-\sin
(x)\,ci(-2n-1/m+1,\zeta x)\right]  ,\\
\int_{0}^{\infty}\frac{\cos(\zeta t)}{(t+x)^{2n+1/m}}dt  &  =\zeta
^{2n+1/m-1}\left[  \cos(\zeta x)\,ci(-2n-1/m+1,\zeta x)+\sin(\zeta
x)\,si(-2n-1/m+1,\zeta x)\right]  \,,
\end{align*}
where%
\begin{align*}
si(\alpha,z)  &  =\int_{z}^{\infty}\frac{\sin(t)}{t^{1-\alpha}}dt,\\
ci(\alpha,z)  &  =\int_{z}^{\infty}\frac{\cos(t)}{t^{1-\alpha}}dt,
\end{align*}
$(\alpha<1)$ the latter integrals being the \textit{generalized sine and
cosine integrals}. \cite{gen} \ Once again we see that the Lommel functions
are related to other special functions.

\end{document}